\documentclass[12pt]{amsart}
\usepackage{amsmath,amssymb}
\usepackage[all]{xy}
\usepackage{pdfpages}
\usepackage[margin=1in]{geometry}
\usepackage{tikz}
\usepackage{tikz-3dplot}
\usetikzlibrary{decorations.pathreplacing, calligraphy}
\usepackage{scalerel}
\usepackage{hyperref}
\usepackage{tkz-euclide}
\usepackage{svg}
\usepackage{multirow}
\usepackage{caption}
\usepackage{subcaption}

\theoremstyle{plain}
\newtheorem{propn}{Proposition}[section]

\newtheorem{example}[propn]{Example}

\theoremstyle{definition}
\newtheorem*{maint}{Main Theorem}
\newtheorem{defn}[propn]{Definition}

\newtheorem*{rem}{Remark}

\theoremstyle{remark}

\DeclareMathOperator{\del}{\partial}

\usepackage{hyperref} 
\hypersetup{pdfborder={0 0 0}, colorlinks=true, urlcolor=blue, citecolor=black}

\def\sideremark#1{\ifvmode\leavevmode\fi\vadjust{\vbox to0pt{\vss
\hbox to 0pt{\hskip\hsize\hskip1em%
\vbox{\hsize2cm\tiny\raggedright\pretolerance10000%
\noindent {\color{red}{#1}}\hfill}\hss}\vbox to8pt{\vfil}\vss}}}%

\numberwithin{equation}{section}

\title{Face embeddings of Archimedean Solids}
\author{Tommy Murphy}
\address[T. Murphy]{Department of Mathematics, California State University Fullerton, 800 N. State College Blvd., Fullerton 92831 CA, USA}
\email{tmurphy@fullerton.edu}
\urladdr{www.fullerton.edu/math/faculty/tmurphy/}
\author{David Weed}
\address[D. Weed]{Department of Mathematics, University of California Riverside}
\email{david@davidweed.net}
\urladdr{www.davidweed.net}

\begin{document}
\begin{abstract}
We characterize the Archimedean solids among uniform polyhedra via face-embeddings into a regular Tetrahedron. This result has been listed without proof in the literature.  
\end{abstract}

\maketitle

To understand any class of objects in mathematics one must describe how one object of the class ``sits inside" another. Define a polyhedron as the convex hull of a finite set of points in $\mathbb{R}^3$ (the vertices). By  definition, all polyhedra in this work are convex and we will interchangeably refer to them as solids. We are concerned with the class of (convex) uniform polyhedra possessing particular  symmetries. The simplest example of a polyhedron is the Tetrahedron ($T$) and our primary goal is to describe how the Archimedean solids can be placed inside $T$ in a natural geometric manner.

 To set the scene, we first recall some  definitions. See \cite{coxeter, coxeter2, crom, grun, hartshorne} for further details and motivation. The main features of a polyhedron are its vertices, edges, and faces which are of dimension 0, 1, and 2 respectively. Recall a polygon is regular if it has all equal side lengths and equal internal angles. A face of a polyhedron is regular if it is a regular polygon, and a  polyhedron is said to be vertex-transitive if every vertex can be mapped to any other vertex isometrically under a global symmetry of the solid.

\begin{defn} A polyhedron $S$ is said to be {\textit{uniform}} if it has regular faces and is vertex-transitive.
\end{defn}

This class includes the regular polyhedra,  also known as the Platonic solids. There are exactly 5 regular polyhedra: the Tetrahedron ($T$), Cube ($C$), Octahedron ($O$), Dodecahedron ($D$), and Icosahedron ($I$). The wider class of uniform polyhedra contains both the Platonic and Archimedean solids plus a few more. As is well-known \cite{grun}, it is a subtle and technical question to precisely define the Archimedean solids. There are issues with chirality, clashes between local and global conditions, and idiosyncratic conventions by different authors  all contributing to  confusion in the literature. The usual, somewhat clunky, definition is that the Archimedean solids are the uniform polyhedra  which are left over when we exclude the standard examples (Platonic solids, prisms and anti-prisms). In this work our goal is to give a geometric characterization of Archimedean solids inside the class of uniform polyhedra.

To this end, let us begin with the main definition in this paper. Given a solid $S$, denote by $\del S$ the boundary of $S$ (i.e. the union of its faces).  
\begin{defn} Let $S_1$ and $S_2$ be two distinct uniform polyhedra. Then $S_1$ admits a ``\textit{k-face embedding}" in $S_2$, written  $S_1\subset_k S_2$, if  $S_1$ is a proper subset of $S_2$  and  $\del S_1 \cap \del S_2$ consists of precisely $k$ faces of $S_1$. 
\end{defn}
 
\begin{maint}
     Let $S$ be a  uniform polyhedron. Then $S\subset_4T$ if, and only if, $S$ is the Icosahedron ($I$), the Octahedron ($O$), or an Archimedean solid. 
\end{maint}

Our interest in this condition was piqued by reading \cite{pugh}, which lists this condition as a characterization of Archimedean solids.   No proof is provided, but Pugh is adamant that this is a fundamental fact, going as far as to use this property to actually \textit{define} the Archimedean solids, presumably as a more appealing alternative to the more technical standard definition we have already discussed.   At various other places in the literature, this result is mentioned \cite{grun, kappraff} but no details have been provided. 

The proof is more important than the theorem. Classically the Archimedean solids are constructed by starting with the Platonic solids and applying various truncations and strictly more complicated operations such as  snubification and  expansion. We refer the reader to \cite{conway, kepler, pugh, stott}.  Moreover, topological arguments are employed as part of the classifications because these operations result in faces where the edge lengths can vary. A cumbersome rescaling of certain faces is then required to produce a uniform solid (c.f. the truncated Cuboctahedron and truncated Icosidodecahedron).  This issue was apparent even to Kepler \cite{kepler}. Kappraff \cite{kappraff} mentions that it is possible to construct all the Archimedean solids in an entirely metric fashion by slicing the Platonic solids with judiciously chosen planes, but does not give precise details. We explicitly describe how to do this and thus eliminate these ambiguous rescaling issues.  Throughout we will use Conway's polyhedra notation \cite{hart}. The notation consists of lowercase letters representing the operation and uppercase letters for the solid on which we perform the operation e.g. $tC$ is the truncated Cube. Sometimes there is a difference  between the operation denoted and the one performed in our construction,  however each solid will be named prior to its notation.

To establish the main theorem we need to explicitly construct a 4-face embedding of each Archimedean solid.  The starting point is the classically known fact that  $O$ and  $I$ admit a 4-face embedding into $T$ \cite{kappraff}. We will build almost all the Archimedean solids from this embedding of $O$ or $I$ using various truncations that do not affect the four faces coincident with the faces of the Tetrahedron. The only exception to this approach is tT, the truncated Tetrahedron, where the 4-face embedding is immediate.
 
\textit{Acknowledgements: We thank CSUF for supporting undergraduate research and supporting D.W. with an MSRP grant. We especially thank the referee, whose insightful comments greatly improved this paper.}

\bigskip

\begin{centering}

\begin{figure}[h]
\includegraphics[width=\linewidth]{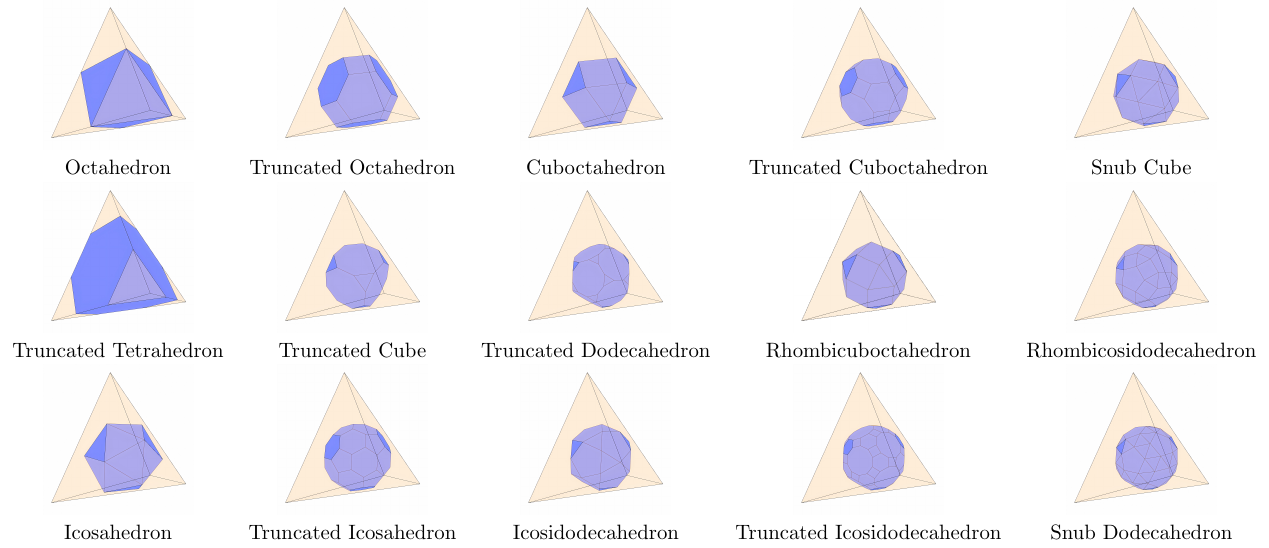}
\caption{Uniform Solids Embedded in a Tetrahedron}
\end{figure}
\end{centering}


\section{Truncations} 

Throughout this work we will start with a uniform polyhedron and apply truncations. We sometimes call our initial polyhedron the ``seed" solid.
\subsection{Vertex Truncation} Take a vertex of a  uniform polyhedron.  We will say that we truncate this vertex if we slice along the solid with a plane orthogonal to the position vector of the vertex. It follows from the classification of uniform polyhedra \cite{coxeter2} that  all vertices of a uniform polyhedron lie on a common circumscribing sphere centered at the origin. Hence vertex truncation  is well defined.  The new face resulting from this procedure will be a regular polygon.  This follows from standard considerations using so-called vertex figures and a straightforward argument using similar triangles (see  \cite{coxeter},  p. 16  for full details).

As vertex truncation is obtained by cutting along planes perpendicular to a vertex,  the operation is generally measured by a length along the edge of the solid rather than the actual depth of cut. Without loss of generality, we scale $O$ and $I$ so that the edges have length $1$. This is convenient as when we truncate $\alpha$ of the edge length, we can simply say we \textit{vertex truncate by $\alpha$}, where $\alpha \in (0,1)$. Vertex truncating to the halfway point of the edge (vertex truncating by $1/2$) is also known as \textit{rectification}.

\begin{example} Starting from one of its vertices of $T$, move along the three edges attached by the truncation value. The points on the edges are coplanar and lay on the vertices of an equilateral triangle. More specifically, the plane defined by these points is parallel to the plane tangent to the circumsphere at the original vertex.
\end{example}

As is well known, starting with the Platonic solids one may vertex truncate by $1/3$ to produce three of the Archimedean solids, namely the truncated Tetrahedron ($tT$), truncated Octahedron ($tO$), and truncated Icosahedron ($tI$). A classical fact is that truncation by $1/2$ (rectification) produces the Icosidodecahedron ($ID$) and the Cuboctahedron ($CO$).

The standard construction of the truncated Cube ($tC$) and truncated Dodecahedron ($tD$) is to truncate $C$ and $D$, respectively, by $1/3$. However, in order to preserve the natural 4-face embedding for these two solids  we  will use an alternate construction. Due to the dual nature of $O$ and $C$ (as well as $I$ and $D$) the cutting planes for the vertex truncation of $O$ are the faces of $C$ and vice-versa (similarly for $I$ and $D$). This means we can reach these required solids by truncating past $1/2$ for $O$ and $I$. Explicitly, vertex truncating $O$ by $\frac{2+\sqrt{2}}{3+2\sqrt{2}}$ produces $tC$ and vertex truncating $I$ by $\frac{2 + \phi}{3 + 2\phi}$ where, $\phi = \frac{1+\sqrt5}{2}$ (the golden ratio), produces $tD$. These values are notable since $\sqrt2$ and $\phi$ are the lengths of the face diagonals of $C$ and $D$.

\subsection{Edge Truncation} Analogous to vertex truncation, edge truncation is defined by cutting by a plane parallel to the tangent plane of the mid-sphere (i.e. the sphere touching each edge midpoint). We measure perpendicular to the edge along the face to determine the depth of cut, so we will informally say we edge-truncate by $\beta$.

The next solids to construct are the Rhombicuboctahedron ($eO$), Rhombicosidodecahedron ($eI$), truncated Cuboctahedron ($tCO$), and truncated Icosidodecahedron ($tID$). To achieve this we perform vertex truncation and edge truncation simultaneously. The notations $tCO$ and $tID$ are not wholly accurate since they cannot be produced by simple vertex truncation as rescaling certain edges would be required. The notation $eO, eI$ comes from the fact that these polyhedra are typically constructed via expansion.  In fact, Boole-Stott \cite{stott} originally produced 11 of the 13 Archimedean solids through vertex, edge, and face expansion. For instance, her construction yields $tCO$ without the issue of having to rescaling certain edges to ensure a uniform solid. However, the expansion operation does not as naturally preserve the 4-face embedding. 

Performing the following vertex and edge truncations of $O$ and $I$ produces $eO, eI, tCO$, and $tID$ while preserving the 4-face embedding. As before, $\phi = \frac{1+\sqrt5}{2}$ denotes the golden ratio.


\begin{center}
\begin{table}[h]
\begin{tabular}{|c|c|c|c|}\hline
   Seed & Solid &  Vertex & Edge \\\hline\hline
   \multirow{2}{*}{$O$} & $eO$  &  $\frac{2}{3+\sqrt{2}_{\mathstrut}}$   & $\frac{\sqrt{3}^{\mathstrut}}{6+2\sqrt{2}}$ \\\cline{2-4}
         & $tCO$ &  $\frac{2+\sqrt{2}}{3+3\sqrt{2}_{\mathstrut}}$ & $\frac{\sqrt{3}^{\mathstrut}}{6+6\sqrt{2}}$ \\\hline
   \multirow{2}{*}{$I$} & $eI$  &  $\frac{2}{3+\phi_{\mathstrut}}$    & $\frac{\sqrt{3}^{\mathstrut}}{6+2\phi}$ \\\cline{2-4}
         & $tID$ &  $\frac{2+\phi}{3+3\phi_{\mathstrut}}$   & $\frac{\sqrt{3}^{\mathstrut}}{6+6\phi}$ \\\hline
\end{tabular}
\caption{Parameters for simultaneous vertex and edge truncation.}
\label{tableve}
\end{table}
\end{center}

\section{Explicit Computations}

Here we present the details of the construction of $tCO$ by applying vertex and edge truncations to $O$ using the values given in the table. For the remaining solids mentioned, the proof is entirely analogous, and details are left to the reader. Starting with $O\subset_4 T$, the first step is to mark the cutting lines on the triangular faces. The red lines represent the planes for vertex truncation by $r_1$ and the green lines denotes edge truncation by $r_2$. In the center the shaded hexagon is a new face of $tCO$. The four faces of $O$ which are subsets of the faces of the circumscribing $T$ will give us the 4-face embedding desired.

\begin{center}
\begin{tikzpicture}[scale=1.3]
\clip(-0.2,-0.2) rectangle (4.2, 3.66);
\draw [line width=1pt] (0,0) -- (4,0) -- (2,3.4641016151377544) -- cycle;
\draw [fill=gray, opacity=0.3] (1.6,1.83) -- (2.4,1.83) -- (2.79,1.16) -- (2.39,0.47) -- (1.61,0.47) -- (1.21,1.16) -- cycle;
\draw [color=green] (1.73, 2.99) -- (3.46, 0);
\draw [color=green] (2.27, 2.99) -- (0.54, 0);
\draw [color=green] (0.27,0.47) -- (3.73, 0.47);
\draw [color=red] (1.06, 1.83) -- (2.94, 1.83);
\draw [color=red] (3.06, 1.63) -- (2.12, 0);
\draw [color=red] (1.88,0) -- (0.94, 1.63);
\draw [decorate, decoration = {brace, raise=1pt}] (1.6, 1.83) -- node [above] {$l_1$} (2.4,1.83);
\draw [decorate, decoration = {brace, raise=1pt, mirror}] (2.4,1.83) -- node [below left] {$l_2$} (2.79,1.16);
\draw [decorate, decoration = {brace, raise=2pt}] (2,3.46) -- node [above right] {$r_1$} (2.94,1.83);
\draw [decorate, decoration = {brace, raise=1pt}] (2.65,1.39) -- node [above=2pt, xshift=-1ex, yshift=-1pt] {$r_2$} (3.06,1.63);
\draw [dashed] (2.65,1.39) -- (3.06,1.63);
\end{tikzpicture}
\end{center}

Firstly we parameterize $l_1$ and $l_2$ in terms of $r_1$ and $r_2$. We obtain $l_1$ by observing the top triangle with side labeled $r_1$ is equilateral. 

\begin{center}
\begin{tikzpicture}[scale=1.5]
\clip(-0.2,-0.2) rectangle (4.2, 3.66);
\draw [line width=1pt] (0,0) -- (4,0) -- (2,3.4641016151377544) -- cycle;
\draw [color=green] (1.73, 2.99) -- (3.46, 0);
\draw [color=green] (2.27, 2.99) -- (0.54, 0);
\draw [color=green] (0.27,0.47) -- (3.73, 0.47);
\draw [color=red] (1.06, 1.83) -- (2.94, 1.83);
\draw [color=red] (3.06, 1.63) -- (2.12, 0);
\draw [color=red] (1.88,0) -- (0.94, 1.63);
\draw [decorate, decoration = {brace, raise=1pt}] (1.6, 1.83) -- node [above] {$l_1$} (2.4,1.83);
\draw [decorate, decoration = {brace, raise=2pt}] (2,3.46) -- node [above right] {$r_1$} (2.94,1.83);
\draw [decorate, decoration = {brace, raise=1pt}] (2.65,1.39) -- node [above=2pt, xshift=-1ex] {$r_2$} (3.06,1.63);
\draw [densely dashed] (2.65,1.39) -- (3.06,1.63);
\end{tikzpicture}
\begin{tikzpicture}[scale=1.5]
\clip(0.86,-0.2) rectangle (3.14, 3.66);
\draw [line width=1pt] (1.06, 1.83) -- (2,3.46) -- (2.94, 1.83);
\draw [color=red] (1.06, 1.83) -- (2.94, 1.83);
\draw [line width=1pt, color=blue] (2.94, 1.83) -- (2.53, 1.6) -- (2.4, 1.83) -- cycle;
\draw [color=green] (1.73, 2.99) -- (2.4,1.83);
\draw [color=green] (2.27, 2.99) -- (1.6, 1.83);
\draw [decorate, decoration = {brace, raise=1pt, mirror}] (1.6, 1.83) -- node [below, yshift=-2pt] {$l_1$} (2.4,1.83);
\draw [decorate, decoration = {brace, raise=2pt}] (2,3.46) -- node [above right] {$r_1$} (2.94,1.83);
\draw [decorate, decoration = {brace, raise=1pt, mirror}] (2.53, 1.6) -- node [below right] {$r_2$} (2.94, 1.83);
\end{tikzpicture}\hspace{-13ex}
\begin{tikzpicture}[scale=1.5]
\draw [line width=1pt, color=blue] (0,0) -- (0.27,0) -- (0, 0.47) -- cycle;
\draw [decorate, decoration = {brace, raise=2pt}] (0,0) -- node [left, xshift=-2pt] {$r_2$} (0, 0.47);
\draw [decorate, decoration = {brace, raise=2pt, mirror}] (0.27,0) -- node [above right] {$\frac{2}{\sqrt{3}}r_2$} (0, 0.47);
\draw [decorate, decoration = {brace, raise=2pt, mirror}] (0,0) -- node [below, yshift=-2pt] {$\frac{r_2}{\sqrt3}$} (0.27,0);
\end{tikzpicture}
\end{center}
Equating sides of the upper triangle yields $r_1 = l_1 + 2\left(\frac{2}{\sqrt{3}}r_2\right)$, or equivalently 
\begin{equation}\label{e1}
l_1 = r_1 -\frac{4}{\sqrt{3}}r_2
\end{equation}
  
To determine $l_2$ we focus on the following blue trapezoid: 

\begin{center}
\begin{tikzpicture}[scale=1.5]
\clip(-0.2,-0.2) rectangle (4.2, 3.66);
\draw [line width=1pt] (0,0) -- (4,0) -- (2,3.4641016151377544) -- cycle;
\draw [color=green] (1.73, 2.99) -- (3.46, 0);
\draw [color=green] (2.27, 2.99) -- (0.54, 0);
\draw [color=green] (0.27,0.47) -- (3.73, 0.47);
\draw [color=red] (1.06, 1.83) -- (2.94, 1.83);
\draw [color=red] (3.06, 1.63) -- (2.12, 0);
\draw [color=red] (1.88,0) -- (0.94, 1.63);
\draw [line width=1pt, color=blue] (2.4, 1.83) -- (2.94, 1.83) --  (3.06, 1.63) -- (2.79, 1.16) -- cycle;
\draw [decorate, decoration = {brace, raise=2pt}] (2,3.46) -- node [above right] {$r_1$} (2.94,1.83);
\draw [decorate, decoration = {brace, raise=1pt}] (2.65,1.39) -- node [above=1.5pt, xshift=-1ex] {$r_2$} (3.06,1.63);
\draw [decorate, decoration = {brace, raise=1pt, mirror}] (2.4,1.83) -- node [below left] {$l_2$} (2.79,1.16);
\draw [densely dashed] (2.65,1.39) -- (3.06,1.63);
\end{tikzpicture}
\begin{tikzpicture}[scale=3]
\draw [line width=1pt, color=blue] (0,0) -- (0.78,0) -- (0.51, 0.47) -- (0.27,0.47) -- cycle;
\draw [opacity=0.4] (0.27,0) -- (0.27,0.47);
\draw [opacity=0.4] (0.51,0) -- (0.51,0.47);
\draw [decorate, decoration = {brace, raise=1pt, mirror}] (0,0) -- node [below] {$l_2$} (0.78,0);
\draw [decorate, decoration = {brace, raise=1pt, mirror}] (0.51, 0.47) -- node [above] {\footnotesize$1-2r_1$} (0.27,0.47);
\draw [decorate, decoration = {brace, raise=1pt, mirror}] (0.78,0) -- node [right] {\footnotesize$\frac{2}{\sqrt{3}}r_2$} (0.51, 0.47);
\draw [decorate, decoration = {brace, raise=5pt, mirror}] (0,0) -- node [below=5pt] {\footnotesize$\frac{r_2}{\sqrt3}$} (0.27,0);
\draw [decorate, decoration = {brace, raise=5pt, mirror}] (0.51,0) -- node [below=5pt] {\footnotesize$\frac{r_2}{\sqrt3}$} (0.78,0);
\end{tikzpicture}
\end{center}
The top value of $1-2r_1$ comes from the fact that $O$ has edge lengths of 1 and we vertex truncate by $r_1$. Clearly
\begin{equation}\label{e2}
l_2 = 1-2r_1 + \frac{2}{\sqrt{3}}r_2.
\end{equation}

On the other hand, $l_2$ is an edge of our new solid and  is also be the edge of the adjacent rectangular face. This new rectangular face is going to be coincident to the edge truncation plane and we can determine its other side length $l_3$ using the following diagram. 
\begin{center}
\begin{tikzpicture}[scale=1.5]
\clip(-0.2,-0.4) rectangle (4.2, 1.61);
\draw [line width=1pt] (0,0) -- (4,0) -- (2,1.41) -- cycle;
\draw [decorate, decoration = {brace, raise=2pt}] (0,0) -- node [above left] {$r_2$} (2,1.41);
\draw [decorate, decoration = {brace, raise=2pt, mirror}] (4,0) -- node [above right] {$r_2$} (2,1.41);
\draw (2, 1.41) node [below] {$\theta$};
\draw [decorate, decoration = {brace, raise=2pt, mirror}] (0,0) -- node [below, yshift=-2pt] {$l_3$} (4,0);
\end{tikzpicture}
\end{center}

This diagram is a cross section of $O$ perpendicular to the edge, where the top vertex would be the midpoint of the edge. To determine the length $l_3$ of the new edge created after edge truncating by $r_2$, an easy application of the Law of Cosines yields
\begin{equation}
    \label{locos}
    l_3 = r_2\sqrt{2}\sqrt{1 - \cos\theta}
\end{equation}
where $\theta$ is the dihedral angle of the seed solid. In our present setting the solid is $O$, so the dihedral angle is $\arccos{(-1/3)}$. Substituting this in yields,  
\begin{equation} \label{e3}
l_3 = \frac{2\sqrt{2}}{\sqrt{3}}r_2.
\end{equation}

Equations (\ref{e1}), (\ref{e2}) and (\ref{e3}) express $l_1$, $l_2$ and $l_3$ as functions of $r_1$ and $r_2$. Now to construct $tCO$, the requirement that the hexagon is regular forces $l_1=l_2$, which from Equations (\ref{e1}) and (\ref{e2}) results in  
$$r_1 = \frac{1}{3} + \frac{2}{\sqrt{3}}r_2.$$
Similarly, looking at the rectangular  face adjacent to the top edge of the hexagon we must have $l_1 = l_3$ to ensure that the resulting solid is  Archimedean. From Equations (\ref{e1}) and (\ref{e3}) this forces
$1-2r_1 + \frac{2r_2}{\sqrt{3}} = \frac{2\sqrt{2}}{\sqrt{3}}r_2$. Solving this system of equations results in $r_1 = \frac{2 + \sqrt{2}}{3 + 3\sqrt{2}}$ and $r_2= \frac{\sqrt{3}}{6 + 6\sqrt{2}}$. 

As already mentioned these simultaneous truncations also produce a rectangular edge with side lengths $l_2$ and $l_3$. A direct computation, which we leave to the reader, shows that simultaneously vertex and edge truncating by the given values for $r_1$ and $r_2$ forces $l_2=l_3,$ and we conclude we have consistency and thus  have constructed $tCO$ as claimed.

\subsection{Skew Truncation}
The final two Archimedean Solids to construct are the snub Cube ($sC$) and the snub Dodecahedron ($sD$). Either of the chiral forms can be produced by the following approach. The desired embedding follows from a construction due to Rotg\'{e} \cite{rotge} (see also \cite{kappraff}) which we call skew truncation.  Starting with $O$ or $I$ we subdivide each edge into two parts that have ratios as per Table \ref{tableve}. The point where the edge is divided is connected to the opposite vertex. This will yield three lines on the face of the original solid intersecting at three points of a regular triangle that has been rotated on the face. Our cutting plane is defined by selecting one edge of this new triangle and the point of the new triangle opposite the original edge. This describes a unique plane for each edge of the original solid. We also define a cutting plane by selecting the vertex of the rotated triangle and the vertices of the rotated triangles around the adjacent original vertex. The key point is that this operation preserves a section of the faces of the original solids, namely the rotated triangle. This means our 4-face embedding is preserved throughout this operation.
\bigskip
\begin{center}
\begin{tabular}{|c|c|c|}\hline
   Solid & Exact Ratio & Approximate \\\hline\hline
   $sC ^{\mathstrut}_{\mathstrut}$& $\frac{1}{3}(1 + \sqrt[3]{19 - 3 \sqrt{33}}^{\mathstrut} + \sqrt[3]{19 + 3 \sqrt{33}})$ & 1.839286755... \\\hline
   $sD ^{\mathstrut}_{\mathstrut}$& $\frac{1}{3} + \frac{2^{5/3} (1 + i \sqrt{3}))}{3 \sqrt[3]{-49 - 27 \sqrt{5} + 3 \sqrt{3 (186 + 98 \sqrt{5})}}}  + \frac{(1 - i \sqrt{3}) \sqrt[3]{-49 - 27 \sqrt{5} + 3 \sqrt{3 (186 + 98 \sqrt{5})}}_{\mathstrut}^{\mathstrut}}{6\cdot 2^{2/3}}$ & 1.943151259...\\\hline   
\end{tabular}
\label{table2}
\end{center}
\bigskip

The approximate values here are pulled from Rotg\'{e} and the exact values are the results of solving the following cubic,
$$r^3-r^2-r-1+2\cos\alpha = 0$$
for $\alpha = 90^\circ,\ 108^\circ$ as shown in \cite{rotge}. For this formula we change to measuring angles in degrees to align with Rotg\'e's conventions. 
The value provided for $sC$ is known as the Tribonacci Constant whereas  the value for $sD$ does not appear to have an analogous title.

\section{Proof of the Main Theorem}

As we have explicitly established the 4-face embeddings required, it remains to show the converse statement. Specifically, we claim that the remaining  uniform polyhedra are not 4-face embeddable in $T$. The other members of this class are the remaining Platonic solids $C$ and $D$, the prisms $P_n$, and the antiprisms $A_n$. Our argument uses the dihedral angles produced by various faces. Observe that when a solid $S \subset_4 T$, there must be four faces of $S$ that have pairwise dihedral angles equal to $\arccos{(1/3)}$, the dihedral angle of $T$. This condition was satisfied in each of the previous constructions because we started with $T$ as our base solid.

For each of these solids, the possible angles between pairs of faces are as follows.

\begin{center}
$\begin{array}{|c|l|} \hline
    \text{Solid}        ^{\mathstrut}_{\mathstrut}& \text{Possible Face Angles} \\\hline\hline
    \text{Cube}         ^{\mathstrut}_{\mathstrut}& \pi/2,\  \pi\\\hline
    \text{Dodecahedron} ^{\mathstrut}_{\mathstrut}& \arccos{(\pm\sqrt{5}/5)},\  \pi\\\hline
    \text{n-Prism}      ^{\mathstrut}_{\mathstrut}& (n-2)k\pi/n: 1 \leq k \leq \lfloor n/2\rfloor,\  \pi/2,\  \pi\\\hline
    \text{n-Antiprism}  ^{\mathstrut}_{\mathstrut}& \arccos(\pm 1/\sqrt{3}\tan(\pi/2n))\\\hline
\end{array}$
\end{center}
This table shows that none of the attainable angles are equivalent to $\arccos{(1/3)}$. This is immediate for all the families except for the antiprisms. An antiprism $A_n$ has two $n$-gonal faces and $2n$ triangular faces. The two possible angles are between the $n$-gonal face and the edge (or vertex) adjacent triangular face. The only way that the angles coincide is when $n=3.$ Here, $A_3$ is equivalent to $O$ and we have already covered this case. The only final possibility is if $T$ is not coincident with the $n$-gonal face and only touches the triangular faces. This implies opposite triangular faces of $A_n$ coincide with $T$ and thus $2\arccos(\pm 1/\sqrt{3}\tan(\pi/2n)) = \arccos{(1/3)}$ for some $n$, which is never true.

The main theorem now follows.  \qed

\begin{rem} Analyzing the face vectors (i.e. the outward normals) was essential in the  proof of the Main Theorem. A necessary condition for a 4-face embedding  of $S_1$ inside $S_2$ to exist is that there are four face vectors of $S_1$ and corresponding face vectors of  $S_2$ which agree. To prove one direction of the Main Theorem,  we have just shown that this is not the case for the appropriate Platonic solids as well as the prisms and anti-prisms.  It is therefore natural to wonder if one could prove our main theorem  purely by analyzing the face vectors of each given Archimedean solid. In our view such an approach would not be as aesthetically satisfying as the constructive proof we provided. Moreover, there {\textit{are}} Archimedean solids with this property which do not admit a 4-face embedding in either direction.  An explicit counterexample is given by taking the truncated Icosahedron  $tI\subset_4 T$ (blue in the picture below) and the truncated Cuboctahedron $tCO\subset_4 T$. Both these solids share four face vectors (corresponding to their respective 4-face embedding into $T$).  It is clear neither solid is contained in the other from the picture.  We conclude that knowing four  face vectors of each polyhedron agree is not enough to conclude one embeds inside the other. This gives some additional context for why our proof avoids analyzing face vectors alone, as a further case-by-case analysis of each face of the solid would be required. 

\begin{centering}
\begin{figure}[h]
\includegraphics[width=.4\linewidth]{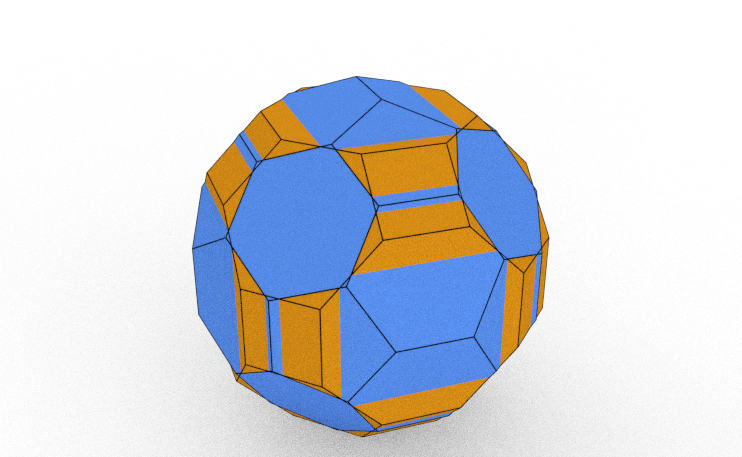}
\caption{$tI$ and $tCO$ 4-face embedded in the tetrahedron $T$.}
\end{figure}
\end{centering}
 
 \end{rem}
\bibliographystyle{plain}
\bibliography{refs}
\end{document}